\documentclass[11pt]{article}
\usepackage{amsfonts}
\usepackage{amssymb}
\usepackage{amssymb,amsfonts,amsmath,amsthm,cite, color}
\usepackage{epsfig}
\newtheorem{theorem}{Theorem}[section]

\newtheorem{e-proposition}[theorem]{Proposition}

\newtheorem{e-definition}[theorem]{Definition\rm}

\newcommand{\Z}{\mathbb{Z}}

\def\H{\mathbb{H}}

\def\og{\leavevmode\raise.3ex\hbox{$\scriptscriptstyle\langle\!\langle$~}}
\def\fg{\leavevmode\raise.3ex\hbox{~$\!\scriptscriptstyle\,\rangle\!\rangle$}}

\begin{document}

\begin{center}
{{\Large\bf Two congruences involving Andrews-Paule's broken 3-diamond partitions and 5-diamond partitions}}

\vskip 6mm

{ Xinhua Xiong
\\[%
2mm] Department of Mathematics\\
China Three Gorges University, Yichang, Hubei Province  443002,
P.R. China \\[3mm]
XinhuaXiong@ctgu.edu.cn \\[0pt%
] }
\end{center}
\begin{abstract}

In this note, we will give proofs of two congruences involving broken 3-diamond partitions and broken
5-diamond partitions which were conjectured by Peter Paule and Silviu Radu.

\end{abstract}

\section{Introduction}
In 2007 Gorges E. Andrews and Peter Paule \cite{And1} introduced a new class of combinatorial objects called
broken k-diamonds. Let $\Delta_k(n)$ denote the number of broken k-diamond partitions of $n$, then they showed that
$$
\sum_{n=0}^{\infty}\Delta_k(n)q^n=
\prod_{n=1}^{\infty}\frac{(1-q^{2n})(1-q^{(2k+1)n})}{(1-q^n)^3(1-q^{(4k+2)n})}.
$$
In 2008 Song Heng Chan \cite{Chan} proved an infinite family of congruences when $k=2$, in 2009 Peter Paule and Silviu Radu \cite{Pau}
 gave two non-standard infinite families of broken 2-diamond congruences. Moreover they stated four conjectures
 related to broken 3-diamond partitions and 5-diamond partitions. In this note we show that their first conjecture and the third conjecture are true:
\begin{theorem}[\cite{Pau}, Conjecture 3.1]\label{1.1}
$$
\prod_{n=1}^\infty(1-q^n)^4(1-q^{2n})^6 \equiv  6\sum_{n=0}^\infty
\Delta_3(7n+5)q^n \pmod{7}.
$$
\end{theorem}
\begin{theorem}[\cite{Pau}, Conjecture 3.3]\label{1.2}
$$
E_4(q^2)\prod_{n=1}^\infty(1-q^n)^8(1-q^{2n})^2 \equiv  8\sum_{n=0}^\infty
\Delta_5(11n+6)q^n \pmod{11}.
$$
\end{theorem}
The techniques in \cite{lovejoy} \cite{Ono96} are adapted here to prove Theorem \ref{1.1} and Theorem \ref{1.2}.

\section{Preliminaries}
Let $\mathbb{H}:=\{z\in\mathbb{C}|{\rm Im}(z)>0\}$ denote the upper half of the complex plane, for a positive integer $N$, the congruence subgroup $\Gamma_0(N)$ of $SL_2(\mathbb{Z})$ is defined by
$$
\Gamma_0(N):=\left\{\left(
                 \begin{array}{cc}
                   a & b \\
                   c & d \\
                 \end{array}
               \right)\Big|ad-bc=1,c\equiv 0\ ({\rm mod}\ N)
\right\}.
$$
 $\gamma=\left(
                 \begin{array}{cc}
                   a & b \\
                   c & d \\
                 \end{array}
               \right)\in SL_2(\mathbb{Z})$ acts on the upper half of the complex plane by the linear
fractional transformation
$
\gamma z:=\frac{az+b}{cz+d}.
$
If $f(z)$ is a function on $\H$, which satisfies
$
f(\gamma z)=\chi (d)(cz+d)^kf(z),
$
where $\chi$ is a Dirichlet character modulo $N$, and $f(z)$ is holomorphic on $\H$ and  meromorphic at all the cusps of $\Gamma_0(N)$, then we call $f(z)$
a weakly holomorphic modular form of weight $k$ with respect to $\Gamma_0(N)$ and character $\chi$. Moreover, if $f(z)$ is holomorphic on $\H$ and at all the cusps of $\Gamma_0(N)$, then we call $f(z)$ a holomorphic modular form of weight $k$ with respect to $\Gamma_0(N)$ and character $\chi$. The set of all holomorphic modular forms of weight $k$ with respect to $\Gamma_0(N)$ and character $\chi$ is denoted by
$\mathcal{M}_k(\Gamma_0(N),\chi)$.

Dedekind's eta function is defined by
$
\eta(z):=q^{\frac{1}{24}}\prod_{n=1}^\infty (1-q^n),
$
where $q=e^{2\pi iz}$ and ${\rm Im}(z)>0$.
A function $f(z)$ is called an eta-product if it can
be written in the form of
$
f(z)=\prod_{\delta |N}\eta^{r_{\delta}}(\delta z),
$
where $N$ and $\delta$ are natural numbers and  $r_{\delta}$ is an integer. The following
Proposition 2.1 obtained by Gordon, Hughes \cite{Gordon84} and Newman \cite{Newman} is useful  to verify  whether an eta-product is a weakly holomorphic modular form.

\begin{e-proposition}[\cite{Ono96}, p.174]\label{prop2.1}
If $f(z)=\prod_{\delta|N}\eta^{r_{\delta}}(\delta
z)$ is an eta-product with
$
k:=\frac{1}{2}\sum_{\delta|N}r_{\delta}  \in \Z
$
satisfying the conditions:
$$
\sum_{\delta|N}\delta r_{\delta}\equiv 0 \ ({\rm mod}\ 24),\quad
\sum_{\delta|N}\frac{N}{\delta} r_{\delta}\equiv 0 \ ({\rm mod}\
24),\quad
$$
then $f(z)$ is a weakly holomorphic modular form of weight $k$ with respect to $\Gamma_0(N)$ with the
character $\chi$, here $\chi$ is defined by $\chi(d)=(\frac{(-1)^ks}{d})$ and $s$ is defined by $s:= \prod_{\delta|N}\delta^{r_{\delta}}$.
\end{e-proposition}

The following Proposition 2.2 obtained by Ligozat \cite{Ligozat75} gives the analytic orders of an eta-product at the
cusps of $\Gamma_0(N)$.

\begin{e-proposition}[\cite{Ono96}, p.174, last line]\label{prop2.2}
Let $c,d$ and $N$ be positive integers with $d|N$ and $(c,d)=1$. If
$f(z)$ is an eta-product satisfying the conditions in Proposition
\ref{prop2.1} for $N$, then the order of vanishing of $f(z)$ at the
cusp $\frac{c}{d}$ is
$$\label{formula}
\frac{N}{24}\sum_{\delta
|N}\frac{(d,\delta)^2r_{\delta}}{(d,\frac{N}{d})d\delta}.
$$
\end{e-proposition}

Let $p$ be a prime, and
$
f(q)= \sum_{n\ge n_0}^{\infty}a(n)q^n
$
be a formal power series, we define
$
U_pf(q)=\sum_{pn\ge n_0}a(pn)q^n.
$
\noindent  If $f(z)\in \mathcal{M}_k(\Gamma_0(N),\chi)$,
then $f(z)$ has an expansion at the point $i\infty$ of the form $f(z)=\sum_{n=n_0}^{\infty}a(n)q^n$
where $q=e^{2\pi iz}$ and ${\rm Im}(z)>0$. We call this expansion the Fourier series
of $f(z)$. Moreover we define $U_pf(z)$ to be the result of applying $U_p$
to the Fourier series of $f(z)$. When it acts on spaces of modular forms and $p|N$, we have
$$
U_p: M_k(\Gamma_0(N),\chi) \rightarrow  \,M_k(\Gamma_0(N),\chi)
$$
In \cite{Sturm} Sturm proved the following criterion to determine  whether two modular forms
are congruent, this reduces the proof of a conjectured congruence to a finite calculation.
In order to state his theorem, we introduce the notion of  the $M$-adic order of a
formal power series.
\begin{e-definition}  Let $M$ be a positive integer and $f = \sum_{n \geq N} a(n) q^n$ be a
formal power series in the variable $q$ with rational integer coefficients. The $M$-adic
order of $f$ is defined by
$$
{\rm Ord}_{\sl{M}}(f) = \mbox{inf}\,\{n \mid a(n)\not\equiv 0\,\,\mbox{mod}\,\,M\}
$$
\end{e-definition}

\begin{e-proposition}[Sturm\cite{Sturm}]\label{prop2.3}
Let $M$ be a positive integer and  $f(z),\, g(z)\in M_k(\Gamma_0(N),\chi)\bigcap \Z[[q]]$.
If
\begin{displaymath}
{\rm Ord}_M(f(z)-g(z))\ge 1+ \frac{kN}{12}\prod_{p}(1+\frac{1}{p}),
\end{displaymath}
 where the product is over
all prime divisors $p$ of $N$. Then $f(z)\equiv g(z)\ ({\rm mod\
}M)$.
\end{e-proposition}
\begin{e-proposition}[Theorem 1.67\cite{Ono04}]\label{prop2.4}
$$
E_4(z)=\frac{\eta^{16}(z)}{\eta^8(2z)}+2^8\frac{\eta^{16}(2z)}{\eta^8(z)},
$$
where $E_4(z)$ is the Eisenstein series of weight 4 for the full modular group.
\end{e-proposition}

\section{Proof of Theorem \ref{1.1} }
\noindent{\it Proof.\textbf{--}}
We define an eta-product
$$
F(z):=\frac{\eta(2z)\eta^9(7z)}{\eta^3(z)\eta(14z)},
$$
\noindent setting $N=56$, we find that $F(z)$ satisfies the conditions of Proposition \ref{prop2.1} and $F(z)$ 
is holomorphic at all cusps of $\Gamma_0(56)$ by using Proposition \ref{prop2.2},
so $F(z)$ is in $\mathcal{M}_3(\Gamma_0(56),\chi)$, where $\chi(d)=(\frac{-1}{d})$ is a Dirichlet character modulo $56$. we note that
$$
F(z)= q^2 \prod_{n=1}^{\infty}\frac{(1-q^{2n})(1-q^{7n})^9}{(1-q^n)^3(1-q^{14n})}
$$
and
$$
\sum_{n=0}^{\infty}\Delta_3(n)q^n=
\prod_{n=1}^{\infty}\frac{(1-q^{2n})(1-q^{7n})}{(1-q^n)^3(1-q^{14n})}.
$$
Applying $U_7$ operator on $F(z)$, we find that
\begin{eqnarray}\label{3}
F(z)|U_7 &=& \left(q^2 \prod_{n=1}^{\infty}\frac{(1-q^{2n})(1-q^{7n})^9}{(1-q^n)^3(1-q^{14n})}  \right)\big|U_7
= \left(q^2 \sum_{n=0}^\infty\Delta_3(n)q^n \prod_{n=1}^\infty(1-q^{7n})^8  \right)\big|U_7\cr\cr
&=& \left(\sum_{n\geq 2}^{\infty}\Delta_3(n-2)q^n  \right)\big|U_7 \cdot \prod_{n=1}^\infty(1-q^{n})^8
= \sum_{7n\geq 2}^{\infty}\Delta_3(7n-2)q^n \prod_{n=1}^\infty(1-q^{n})^8\cr\cr
&=& q\sum_{7n\geq 2}^{\infty}\Delta_3(7n-2)q^{n-1} \prod_{n=1}^\infty(1-q^{n})^8\cr\cr
&=& q\sum_{n\geq 0}^{\infty}\Delta_3(7n+5)q^{n} \prod_{n=1}^\infty(1-q^{n})^8.
\end{eqnarray}

We define another eta-product
$$
G(z):=\frac{\eta^6(2z)\eta^2(7z)}{\eta^2(z)}.
$$
By Proposition \ref{prop2.1} and Proposition \ref{prop2.2}, we find that $G$ is also in $\mathcal{M}_3(\Gamma_0(56),\chi)$, where $\chi(d)=(\frac{-1}{d})$ is a Dirichlet character modulo $56$. Moreover, we have
\begin{eqnarray}\label{4}
G(z)&=&\frac{\eta^6(2z)\eta^2(7z)}{\eta^2(z)}=\eta^{12}(z)\eta^6(2z)\frac{\eta^2(7z)}{\eta^{14}(z)}
\equiv  \eta^{12}(z)\eta^6(2z) \pmod{7}\cr\cr
&\equiv & q \prod_{n=1}^\infty(1-q^{n})^{12}(1-q^{2n})^6\pmod{7}.
\end{eqnarray}
Where we used the elementary fact
$$
\frac{\eta^2(7z)}{\eta^{14}(z)}=\prod_{n=1}^\infty\frac{(1-q^{7n})^{2}}{(1-q^{n})^{14}}\equiv 1 \pmod{7}.
$$
We note that our Theorem \ref{1.1} is equivalent to the congruence:
$$
q \prod_{n=1}^\infty(1-q^{n})^{12}(1-q^{2n})^6 \equiv 6 q\sum_{n\geq 0}^{\infty}\Delta_3(7n+5)q^{n} \prod_{n=1}^\infty(1-q^{n})^8 \pmod{7},
$$
i.e.
$$
G(z) \equiv 6 F(z)|U_7 \pmod{7}.
$$
Using Sturm's theorem \ref{prop2.3}, it suffices to verify the congruence above holds for the
first $\frac{3}{12} \cdot  [ SL_2(\Z):\Gamma_0(56)]+1=25$ terms, which is easily completed by using
Mathematica 6.0.
\section{Proof of Theorem \ref{1.2} }
The proof of Theorem \ref{1.2} is similar. The difference is that we need to construct
 two eta-products to represent the left hand side of the equation in Theorem \ref{1.2} up to a factor by using Proposition \ref{prop2.4}.

\noindent{\it Proof.\textbf{--}}
Define
$$
H(z):=\frac{\eta(2z)\eta^{13}(11z)}{\eta^3(z)\eta(22z)},
$$
Setting $N=88$, we find that $H(z)$ satisfies the conditions of Proposition \ref{prop2.1} and $H(z)$ is holomorphic at
 all the cusps of $\Gamma_0(88)$ by Proposition \ref{prop2.2}, so $H(z)$ is in $\mathcal{M}_5(\Gamma_0(88),\chi)$, where $\chi(d)=(\frac{-1}{d})$ is a Dirichlet character modulo $88$. We note that
$$
H(z)= q^5 \prod_{n=1}^{\infty}\frac{(1-q^{2n})(1-q^{11n})^{13}}{(1-q^n)^3(1-q^{22n})}.
$$
and
$$
\sum_{n=0}^{\infty}\Delta_5(n)q^n=
\prod_{n=1}^{\infty}\frac{(1-q^{2n})(1-q^{11n})}{(1-q^n)^3(1-q^{22n})}.
$$
As before, applying $U_{11}$ operator on $H(z)$, we find that
\begin{eqnarray}\label{5}
H(z)|U_{11} &=& \left(q^5 \prod_{n=1}^{\infty}\frac{(1-q^{2n})(1-q^{11n})^{13}}{(1-q^n)^3(1-q^{22n})}  \right)\big|U_{11}
= \left(q^5 \sum_{n=0}^\infty\Delta_5(n)q^n \prod_{n=1}^\infty(1-q^{11n})^{12}  \right)\big|U_{11}\cr\cr
&=& \left(\sum_{n\geq 5}^{\infty}\Delta_5(n-5)q^n  \right)\big|U_{11} \cdot \prod_{n=1}^\infty(1-q^{n})^{12}
= \sum_{11n\geq 5}^{\infty}\Delta_5(11n-5)q^n \prod_{n=1}^\infty(1-q^{n})^{12}\cr\cr
&=& q\sum_{11n\geq 5}^{\infty}\Delta_5(11n-5)q^{n-1} \prod_{n=1}^\infty(1-q^{n})^{12}\cr\cr
&=& q\sum_{n\geq 0}^{\infty}\Delta_5(11n+6)q^{n} \prod_{n=1}^\infty(1-q^{n})^{12}.
\end{eqnarray}
We define another two eta-products by
$$
L_1(z):=\frac{\eta^{18}(2z)\eta^{2}(11z)}{\eta^2(z)\eta^8(4z)}\quad\text{and}\quad L_2(z):=\frac{\eta^{16}(4z)\eta^{2}(11z)}{\eta^6(2z)\eta^2(z)}.
$$
Setting $N=88$, it is easy to verify that both $L_1(z)$ and $L_2(z)$ satisfy the conditions in Proposition \ref{prop2.1} 
and both are holomorphic at all the cusps of $\Gamma_0(88)$ by using Proposition \ref{prop2.2},
hence both $L_1(z)$ and $L_2(z)$ are in $\mathcal{M}_5(\Gamma_0(88),\chi)$, where $\chi(d)=(\frac{-1}{d})$ is a Dirichlet character modulo $88$. So $L(z):= L_1(z) +2^8L_2(z) $ is in $\mathcal{M}_5(\Gamma_0(88),\chi)$. On the other hand,
\begin{eqnarray}\label{6}
L(z)&=&\frac{\eta^{16}(2z)}{\eta^8(4z)}\cdot \frac{\eta^{2}(2z)\eta^{2}(11z)}{\eta^2(z)}
+2^8\frac{\eta^{16}(4z)}{\eta^8(2z)}\cdot  \frac{\eta^{2}(2z)\eta^{2}(11z)}{\eta^2(z)}\cr\cr
&=&E_4(2z)\cdot \frac{\eta^{2}(2z)\eta^{2}(11z)}{\eta^2(z)}\cr\cr
&=& E_4(2z)\cdot \eta^{20}(z)\eta^2(2z)\cdot \frac{\eta^{2}(11z)}{\eta^{22}(z)}
\equiv   E_4(2z)\cdot \eta^{20}(z)\eta^2(2z) \pmod{11}\cr\cr
&=& E_4(q^2)\cdot q\prod_{n=1}^{\infty}(1-q^{2n})^2(1-q^{n})^{20}.
\end{eqnarray}
We note that Theorem \ref{1.2} is equivalent to the following congruence of modular forms by using the expressions (\ref{5}) and
(\ref{6}):
$$
L(z) \equiv 8 H(z)|U_{11} \pmod{11}.
$$
Using Sturm's theorem \ref{prop2.3}, it suffices to verify the congruence above holds for the
first $\frac{5}{12} \cdot  [ SL_2(\Z):\Gamma_0(88)]+1=61$ terms, which is easily completed by using
Mathematica 6.0.\qed


\end{document}